\documentclass[reqno]{amsart}
\usepackage{amsmath, amssymb, amsthm, geometry, enumerate, graphicx, amsfonts, hyperref, mathrsfs}
\hypersetup{colorlinks=true,linkcolor=red, anchorcolor=blue, citecolor=blue, urlcolor=red, filecolor=magenta, pdftoolbar=true}
\theoremstyle{plain}
\newtheorem{thm}{\it Theorem}[section]

\theoremstyle{remark}
\newtheorem{defn}[thm]{Def{}inition}
\newtheorem{rem}[thm]{Remark}
\newtheorem{exa}[thm]{Example}
\numberwithin{equation}{section}
\usepackage{enumerate, float}
\usepackage{relsize}
\allowdisplaybreaks
\begin{document}
	
	\title [Sums of Frames with Several Generators from the  Weyl--Heisenberg Group]{Sums of Frames from the  Weyl--Heisenberg Group and  Applications to  Frame Algorithm}
	
\author[Divya Jindal]{Divya Jindal}
\address{{\bf{Divya Jindal}}, Department of Mathematics,
University of Delhi, Delhi-110007, India.}
\email{divyajindal193@gmail.com}

\author[Jyoti]{Jyoti}
\address{{\bf{Jyoti}}, Department of Mathematics,
University of Delhi, Delhi-110007, India.}
\email{jyoti.sheoran3@gmail.com}
\author[Lalit   Kumar Vashisht]{Lalit  Kumar  Vashisht$^*$}
\address{{\bf{Lalit  Kumar  Vashisht}}, Department of Mathematics,
University of Delhi, Delhi-110007, India.}
\email{lalitkvashisht@gmail.com}

	\begin{abstract}
The relationship between the frame bounds of  frames  for  the space $L^2(\mathbb{R})$ with several  generators from the Weyl-Heisenberg group and the scalars linked to the sum of frames  is   examined in this paper. We give sufficient conditions for the finite sum of frames of the space $L^2(\mathbb{R})$ from the Weyl-Heisenberg group, with explicit frame bounds, in terms of frame bounds and scalars involved in the finite sum of frames,  to be a frame for $L^2(\mathbb{R})$. We show that the sum of  frames from the Weyl-Heisenberg group and its dual frame always constitutes a frame.  We provide sufficient conditions for the sum of images of frames under bounded linear operators acting on $L^2(\mathbb{R})$ in terms of lower bounds of their Hilbert adjoint operator to be a frame.   The finite sum of frames where frames are perturbed by bounded sequences of scalars is also discussed. As an application of the results, we show that the frame bounds of sums of frames can increase the rate of approximation in the frame algorithm.
	\end{abstract}
	
	\renewcommand{\thefootnote}{}
	
	\footnote{
		The research of Divya  Jindal is supported by the Council of Scientific $\&$ Industrial Research (CSIR), India.  Grant No.: 09/045(1680)/2019-EMR-I. Lalit  is  supported by the Faculty Research Programme Grant-IoE, University of Delhi \ (Grant No.: Ref. No./IoE/2023-24/12/FRP). Jyoti was invited by the third author for research collaboration in Department of Mathematics, University of Delhi.  \\
$^*$Corresponding author: Lalit Kumar Vashisht}

	\maketitle
	
	\baselineskip15pt
	\section{Introduction}
The rise of frame based technologies for signal processing, image processing,  approximation of signals in frame algorithms and also  dozens of applied areas  has attracted mathematicians and engineers. Using a fundamental approach to signal decomposition in terms of elementary signals due to Gabor \cite{Gabor}, Duffin and Schaeffer \cite{DS} introduced the concept of frames for the space $L^2(-\gamma, \gamma)$, $\gamma > 0$. Let $\mathcal{H}$ be a separable (real or complex) Hilbert space with inner product $\langle \cdot , \cdot \rangle$. A countable collection of vectors $\{f_k\}_{k \in \mathcal{I}} \subset \mathcal{H}$ is called a \emph{frame} (or \emph{discrete Hilbert frame}) for $\mathcal{H}$ if there exist finite positive scalars $A_o$ and $B_o$ such that
	\begin{align}\label{framedefn}
		A_o\|f\|^2 \leq \sum_{k\in \mathcal{I}} |\langle f , f_k \rangle |^2 \leq B_o\|f\|^2  \ \text{for all} \ f \in \mathcal{H}.
	\end{align}
The scalars $A_o$ and $B_o$ are called \emph{lower frame bound} and \emph{upper frame bound}, respectively, and are not necessary unique. If $A_o=B_o$, then $\{f_k\}_{k \in \mathcal{I}}$ is called a \emph{tight frame} (or $A_o$-tight frame) and  \emph{Parseval frame} if $A_o=B_o=1$. If only the upper inequality holds in \eqref{framedefn}, we say that $\{f_k\}_{k \in \mathcal{I}}$ is a Bessel sequence with Bessel bound $B_o$. The map $S \colon \mathcal{H} \rightarrow \mathcal{H}$, defined by $S \colon f \rightarrow \sum_{k\in \mathcal{I}} \langle f , f_k \rangle f_k$ is called the \emph{frame operator} which is bounded, linear and invertible on $\mathcal{H}$. This gives the \emph{reconstruction formula} for each $f \in \mathcal{H}$: $f= S S^{-1}f = \sum_{k\in \mathcal{I}}\langle S^{-1}f , f_k \rangle f_k$.
	Thus, a frame for $\mathcal{H}$ allows each element in $\mathcal{H}$ to be expressed as a linear combination, not necessarily unique, of the elements of the frame $\{f_k\}_{k \in \mathcal{I}}$. Tight frames provide the following  representation for each $f$ in $\mathcal{H}$:
	\begin{thm}\cite[Corollary 5.1.7]{ole}\label{tightthm}
If $\{f_k\}_{k\in \mathcal{I}}$ is  $A_o$-tight frame for $\mathcal{H}$, then $f= \frac{1}{A_o}\sum\limits_{k\in \mathcal{I}}\langle f, f_k \rangle f_k$ for all $ f \in \mathcal{H}$.
\end{thm}
Duality  of frames also gives a series representation of each vector of the underlying Hilbert space. Let $\{f_k\}_{k\in \mathcal{I}}$ be a frame for $\mathcal{H}$. A frame $\{g_k\}_{k\in \mathcal{I}}$ for $\mathcal{H}$ is called a \emph{dual frame} of $\{f_k\}_{k\in \mathcal{I}}$, if for all $f\in \mathcal{H}$,
$f= \sum_{k\in \mathcal{I}} \langle f, f_k \rangle g_k \quad \text{or} \quad f= \sum_{k\in \mathcal{I}} \langle f, g_k \rangle f_k$.
Recently, there have been significant advancements in frame theory research, and frames have found widespread applications in various fields, such as image and signal processing, quantum measurements, coding and wireless communication, sampling theory, iterated function systems, wavelet theory, among others. These applications have been extensively studied by researchers, as evidenced by several works, including \cite{ACMT, DV1, DV2, JVash,  VD3, Zalik1, Zalik2}.  Christensen \cite{ole}, Gr\"{o}chenig \cite{KG}, Heil \cite{Heil20}, Heil and Walnut \cite{H89},  Krivoshein, Protasov and Skopina \cite{AK} and Young \cite{Young} are good references for basic theory  of frames.
	\subsection{Related work} 	
 In recent years, Gabor frames have been studied extensively due to their potential applications in time-frequency analysis, in particular, to decompose and stable analysis of signals. There is huge literature on Gabor frames and their applications in many areas of physics and applied mathematics, see \cite{GuidoI, FZUS,  KG, GuidoII, Heil20, H89, lalkv, Zalik1, Zalik2} and many references therein. The coherent states associated with a unitary irreducible representation of a Lie group induce Gabor  and wavelet analysis, we refer to  \cite{KG, H89, AP} for technical details. Specifically, Gabor analysis corresponds to the unitary irreducible representation of the Weyl-Heisenberg group. In a recent study by authors of  \cite{SB}, the relationship between Gabor and wavelet analysis was explored using the concept of ``contraction" between the Weyl-Heisenberg group and extended affine group. They established unitary irreducible representations of the Weyl-Heisenberg group as a contraction of representations of the extended affine group. In  \cite{SB}, the authors established a connection between coherent states, tight frames, and resolutions of identity using the technique of contraction \cite{SB}. In this direction, two authors of this paper, in \cite{JVas}, proved  sufficient conditions for the existence of Gabor frames and wavelet frames with several generators for the Weyl-Heisenberg group and the extended affine group, respectively. They also studied matrix-valued nonstationary frames from the Weyl-Heisenberg group \cite{DivII}.

 It is worth recalling that, in general, a finite linear combination of frames for a given space does  not constitute  a frame for that space. In a recent publication \cite{OSCT}, Obeidat, Samarah, Casazza, and Tremain investigated the concept of finite sums of frames in Hilbert spaces. Specifically, they examined the finite sum of frames that incorporate Bessel sequences and frames under bounded linear operators and frame operators on the underlying space. They gave necessary and sufficient conditions on Bessel sequences $\{f_k\}_{k \in \mathcal{I}}$ and $\{g_k\}_{k \in \mathcal{I}}$ and bounded linear operators $L_1$, $L_2$ on a separable Hilbert space $\mathcal{H}$ so that $\{L_1f_k +L_2 g_k\}_{k \in \mathcal{I}}$ is a frame for $\mathcal{H}$. Recently, sums of frames of  wave packet systems were studied in  \cite{JDLG}. Exploring new relationships and constraints on scalar coefficients and frame bounds for finite sums of frames in the underlying space would be of great interest. On the other hand, frame algorithms are used to approximate signals where frame bounds of a given frame play a  significant role. To be exact, the rate of approximation depends on the frame bounds of a given frame, see Section \ref{subsec2}. Motivated by the above work, we study sums of frames from the Weyl-Heisenberg group. Notable contributions  include sufficient conditions for finite sums of frames in terms of given frame bounds and scalars involved in the sum of frames. We show that if we use the frame bounds of sums of frames, then the rate of approximation in the frame algorithm can be increased considerably.
	
	\subsection{Overwiew of the paper}
	The paper is organised as follows: In Section \ref{sec2}, we give some basic notations and provide an overview on Gabor system and Weyl-Heisenberg group.  Section \ref{sec3} has four main results related to frame conditions for sums of frames from the Weyl-Heisenberg group.  Theorem \ref{sfn} provides sufficient conditions with certain frame bounds for a  finite linear combination  of   frames for the Weyl-Heisenberg group to be a frame. Afterwards, Theorem \ref{dual} shows that sum of a frame with its dual frame is always a frame. A sufficient condition for sum of images of two frames under bounded linear operators to be a frame is given in Theorem \ref{operator}. In Theorem \ref{bddseq} relation between two bounded sequences and bounds of two frames can be found which is justified in Example \ref{bddexa}. Section \ref{subsec2} gives applications of frame bounds of sums of frames in the frame algorithm and convergence analysis. In Example \ref{exfal}--Example \ref{bddexaalgo},  we show that the frame bounds associated with sums of frames can decrease the width of the frame which increases  the rate of convergence  in the frame algorithm.

	\section{Preliminaries}\label{sec2}
	In this paper, we use the following symbols: $\mathbb{N}$ for the set of natural numbers, $\mathbb{Z}$ for the set of integers,   $\mathbb{R}$ for the set of real numbers and $\mathbb{C}$ denotes the set of complex numbers. Re$z$ denotes the real part of a complex number $z$. As is standard,  $L^2(\mathbb{R})$ denotes the space of  equivalence classes of square integrable (in the sense of Lebesgue) functions.  $L^2(\mathbb{R})$  is a Hilbert space with respect to the standard inner product  given by
$\langle f , g \rangle = \int_{\mathbb{R}} f(x) \overline{g(x)} dx$, $f$,  $g \in L^2(\mathbb{R})$ which induces the norm:  $\|f\| =\sqrt{ \langle f ,f\rangle}$,  $f \in L^2(\mathbb{R})$. An operator $V$ acting on  $L^2(\mathbb{R})$ is said to be bounded below if there exists a real constant $M_o$ such that $\|Vf\| \geq M_o \|f\|$ for all $f \in L^2(\mathbb{R})$. The Hilbert-adjoint operator of a bounded linear operator $\Theta$ acting on  $L^2(\mathbb{R})$ is denoted by  $\Theta^*$ and defined by the equation  $\langle \Theta f, g\rangle = \langle  f, \Theta^*g\rangle$ for all $f$, $g \in L^2(\mathbb{R})$.
 The Gabor system or time-frequency shift is  based on the translation and modulation operator on $L^2(\mathbb{R})$ which are defined as follows:
	\begin{align*}
		&   \ T_a: f (x) \mapsto f (x-a) \ \quad (\text{Translation by} \ a \in \mathbb{R} );\\
		& E_b: f (x) \mapsto e^{2\pi ib\cdot x} f(x) \quad (\text{Modulation by} \ b \in \mathbb{R}).
	\end{align*}
	A collection of functions of the form $\mathcal{G}(a, b, f): = \{E_{mb}T_{na}f\}_{m, n \in \mathbb{Z}} = \{e^{2 \pi i m b x} f(x -na)\}_{m, n \in \mathbb{Z}}$ is called a \emph{Gabor system} in $L^2(\mathbb{R})$.	
	
	\subsection{The Weyl-Heisenberg group}
	We define the Weyl-Heisenberg group as the outer semidirect product of the groups $\mathbb{R}^2$ and $\mathbb{R}$, denoted by $\mathbb{R}^2 \rtimes_\psi \mathbb{R} = \mathcal{W}$. Here, $\psi:\mathbb{R} \rightarrow Aut(\mathbb{R}^2)$ is a homomorphism given by $\psi_\lambda((x_1,x_2))=(x_1+\lambda x_2,x_2)$ for every $\lambda,x_1,x_2 \in \mathbb{R}$. We denote the group multiplication law by
	\begin{align*}
		((x_1,x_2),\lambda)((y_1,y_2),\mu)=(\psi_\lambda(y_1,y_2)+(x_1,x_2),\lambda+\mu).
	\end{align*}
	Recall that for every $P\in \mathbb{R}^*, Q \in \mathbb{R}$, the unitary irreducible representation, $\sigma^{P,Q}:\mathcal{W} \rightarrow \mathcal{U}(L^2(\mathbb{R}))$ is  given by
	\begin{align*}
		(\sigma^{P,Q}(((u_1,u_2),c))f)(x)= e^{i[P(u_1+xu_2)+Qu_2]}f(c+x).
	\end{align*}
	
	Let $N$ be a fixed natural number, and for $l \in \{1,2,\dots,N\}$ consider  real numbers $q_0^{(l)}$ and $p_0^{(l)}$ such that $|p_0^{(l)}q_0^{(l)}|< 2\pi$. Consider the discrete subset of $\mathcal{W}$ given by
	\begin{align*}
		\mathcal{W}_{q_0^{(l)},p_0^{(l)}}=\bigg\{\bigg(\bigg(\frac{mnp_0^{(l)}q_0^{(l)}}{2},mp_0^{(l)}\bigg),nq_0^{(l)}\bigg)|n,m \in \mathbb{Z}\bigg\}.
	\end{align*}
	For $0\neq \phi_l \in L^2(\mathbb{R})$ where $l \in \{1,2,\dots,N\}$ and $n,m \in \mathbb{Z}$, the sequence $\Big\{	\phi_{m,n,l}^{P,Q}\Big\}_{m, n\in\mathbb{Z} \atop l\in\{1,2,\dots,N\}}$ defined by
	\begin{align*}
		\phi_{m,n,l}^{P,Q}(x)&=\sigma^{P,Q}\bigg(\bigg(\frac{mnp_0^{(l)}q_0^{(l)}}{2},mp_0^{(l)}\bigg),nq_0^{(l)}\bigg)\phi_l(x)\\
		&= e^{i [\frac{P m n p_0^{(l)}q_0^{(l)}}{2}  + Q m p_0^{(l)}]} e^{i P x m p_0^{(l)}} \phi_l\big(x+ n q_0^{(l)}\big),
	\end{align*}
	 is dense in $L^2(\mathbb{R})$, see \cite{SB} for technical details. It is easy to see that the sequence $\Big\{	\phi_{m,n,l}^{P,Q}\Big\}_{m, n\in\mathbb{Z} \atop l\in\{1,2,\dots,N\}}$ is similar to the Gabor structure.

\begin{defn}
		Let $N$ be a fixed natural number. The collection of functions $\big\{\phi_{m,n,l}^{P,Q}\big\}_{m, n\in\mathbb{Z} \atop l\in\{1,2,\dots,N\}}$ in $L^2(\mathbb{R})$ is said to be a \emph{frame} (or \emph{Gabor frame})  for $L^2(\mathbb{R})$ if there exist  $A, B \in (0, \infty)$ such that
		\begin{align*}
			A \|f\|^2 \leq \sum_{l = 1}^{N}\sum_{m, n\in \mathbb{Z}}|\langle f,\phi_{m,n,l}^{P,Q}\rangle|^2 \leq B\|f\|^2 \ \text{for all} \ f \in L^2(\mathbb{R}).
		\end{align*}
	\end{defn}
As in case of standard frames, scalars $A$, $B$ are called \emph{lower frame bound} and \emph{upper frame bound}, respectively,  of the frame  $\big\{\phi_{m,n,l}^{P,Q}\big\}_{m, n\in\mathbb{Z} \atop l\in\{1,2,\dots,N\}}$.
The dual of a frame  $\big\{\phi_{m,n,l}^{P,Q}\big\}_{m, n\in\mathbb{Z} \atop l\in\{1,2,\dots,N\}}$  for  $L^2(\mathbb{R})$ is a  frame $\big\{\psi_{m,n,l}^{P,Q}\big\}_{m, n\in\mathbb{Z} \atop l\in\{1,2,\dots,N\}}$ of  $L^2(\mathbb{R})$ such that each $ f \in L^2(\mathbb{R})$ can be expressed as:
	\begin{align}\label{dualw}
		f= \sum_{l = 1}^{N}\sum_{m, n\in \mathbb{Z}} \langle  f, \phi_{m,n,l}^{P,Q} \rangle \psi_{m,n,l}^{P,Q} \quad \text{or} \quad
		f= \sum_{l = 1}^{N}\sum_{m, n\in \mathbb{Z}} \langle  f, \psi_{m,n,l}^{P,Q} \rangle \phi_{m,n,l}^{P,Q}
	\end{align}
For fundamental results of dual frames with different structure, we refer to \cite{ole, Heil20, Janssen, JVash}.

\section {Sums of Frames from the Weyl--Heisenberg Group}\label{sec3}
We begin with  a sufficient condition, with explicit frame bounds, for a  finite linear combination  of  frames of the space $L^2(\mathbb{R})$  from  the Weyl-Heisenberg group to be a frame for $L^2(\mathbb{R})$.
	\begin{thm}\label{sfn}
		Let $k$ and $N$ be fixed natural numbers, and $\{c_1, \dots, c_k\} \subset \mathbb{C} \setminus \{0\}$.	For each  $i \in \{1,2,\dots,k\}$, let  $\Big\{\phi_{m,n,l}^{(i)^{P,Q}}\Big\}_{m,n \in \mathbb{Z} \atop l \in \{1,2,\dots,N\}}$ be a frame for $L^2(\mathbb{R})$ with frame bounds $A_i$, $B_i$.  Assume that
		\begin{align} \label{sfc}
			|c_j|A_j+ \sum_{i=1 \atop i\neq j}^{k}\bigg|\frac{c_i^2}{c_j}\bigg|A_i > 2 \sqrt{B_j} \Big(\sum_{i=1 \atop i \neq j}^{k}|c_i|\sqrt{B_i}\Big)
		\end{align}
for some $j \in \{1,2,\dots,k\}$.		Then, the finite sum of frames $\left\{ \sum\limits_{i=1}^{k}c_i\phi_{m,n,l}^{(i)^{P,Q}}\right\}_{m,n \in \mathbb{Z} \atop l \in \{1,2,\dots,N\}}$ is a frame for $L^2(\mathbb{R})$ with frame bounds
		\begin{align*}
			\Big(\sum\limits_{i=1}^{k}|c_i|^2 A_i-2\sum_{i=1 \atop i\neq j}^{k}|c_jc_i|\sqrt{B_j B_i}\Big) \ \ \text{and} \ \  k\sum_{i=1}^{k}|c_i|^2B_i.
		\end{align*}
	\end{thm}
	\proof
	For every  $f\in L^2(\mathbb{R})$, we compute
	\begin{align}\label{st1}
		&\sum_{l=1}^N\sum_{m,n \in \mathbb{Z}} \big|\langle f, \sum_{i=1}^{k}c_i\phi_{m,n,l}^{(i)^{P,Q}}\rangle\big|^2\notag\\
		&= \sum_{l=1}^N \sum_{m,n \in \mathbb{Z}} \big|\sum_{i=1}^{k}\langle f,  c_i\phi_{m,n,l}^{(i)^{P,Q}} \rangle \big|^2 \notag\\
		& \geq \sum_{l=1}^N\sum_{m,n \in \mathbb{Z}} \bigg(\big|\langle f, c_j\phi_{m,n,l}^{(j)^{P,Q}} \rangle\big|- \sum_{i=1 \atop i \neq j}^{k}\big|\langle f, c_i\phi_{m,n,l}^{(i)^{P,Q}} \rangle \big|\bigg)^2 \notag\\
		&= \sum_{l=1}^N\sum_{m,n \in \mathbb{Z}} \bigg(\big|\langle f, c_1\phi_{m,n,l}^{(1)^{P,Q}} \rangle \big|^2 + \big|\langle f, c_2\phi_{m,n,l}^{(2)^{P,Q}} \rangle \big|^2 +\dots \big|\langle f, c_k\phi_{m,n,l}^{(k)^{P,Q}} \rangle \big|^2 \notag \\
		&- 2\sum_{i=1 \atop i\neq j}^{k}\big|\langle f, c_j\phi_{m,n,l}^{(j)^{P,Q}} \rangle\big| \cdot \big|\langle f, c_i\phi_{m,n,l}^{(i)^{P,Q}} \rangle \big| + 2\sum_{i=1 \atop i \neq j}^{k}\sum_{s=1 \atop s\neq i,j}^{k}\big|\langle f, c_i\phi_{m,n,l}^{(i)^{P,Q}} \rangle\big| \cdot \big|\langle f, c_s\phi_{m,n,l}^{(s)^{P,Q}} \rangle\big|\bigg)\notag\\
		& \geq \sum_{l=1}^N\sum_{m,n \in \mathbb{Z}} \bigg(\sum_{i=1}^{k} \big|\langle f, c_i\phi_{m,n,l}^{(i)^{P,Q}} \rangle \big|^2 - 2\sum_{i=1 \atop i\neq j }^{k} \big|\langle f, c_j\phi_{m,n,l}^{(j)^{P,Q}} \rangle\big| \cdot \big|\langle f, c_i\phi_{m,n,l}^{(i)^{P,Q}} \rangle\big|\bigg)\notag\\
		& \geq \sum_{l=1}^N\sum_{m,n \in \mathbb{Z}}\sum_{i=1}^{k} \big|\langle f, c_i\phi_{m,n,l}^{(i)^{P,Q}} \rangle\big|^2  - 2\sum_{i=1 \atop i \neq j}^{k} \bigg(\sum_{l=1}^N\sum_{m,n \in \mathbb{Z}} \big|\langle f, c_j\phi_{m,n,l}^{(j)^{P,Q}} \rangle\big|^2\bigg)^{\frac{1}{2}} \bigg(\sum_{l=1}^N\sum_{m,n \in \mathbb{Z}}\big|\langle f, c_i\phi_{m,n,l}^{(i)^{P,Q}} \rangle\big|^2\bigg)^{\frac{1}{2}}\notag\\
		& \geq \sum_{i=1}^{k}|c_i|^2 A_i \|f\|^2 - 2\sum_{i=1 \atop i \neq j}^{k}|c_jc_i|\sqrt{B_j B_i} \|f\|^2\notag\\
		& = \Big(\sum_{i=1}^{k}|c_i|^2 A_i - 2\sum_{i=1 \atop i \neq j}^{k}|c_jc_i|\sqrt{B_jB_i}\Big) \|f\|^2.
	\end{align}
	This gives the lower frame bound. For the upper frame bound, we have
	\begin{align}\label{st2}
		\sum_{l=1}^N\sum_{m,n \in \mathbb{Z}} \Big|\Big\langle f,\sum_{i=1}^{k}c_i\phi_{m,n,l}^{(i)^{P,Q}}\Big\rangle\Big|^2
		&= \sum_{l=1}^N \sum_{m,n \in \mathbb{Z}}\Big|\sum_{i=1}^{k}\langle f, c_i\phi_{m,n,l}^{(i)^{P,Q}} \rangle\Big|^2 \notag\\
		& \leq k \sum_{i=1}^{k}\sum_{l=1}^N \sum_{m,n \in \mathbb{Z}}\big|\langle f,c_i\phi_{m,n,l}^{(i)^{P,Q}}\rangle\big|^2 \notag\\
		& \leq k\sum_{i=1}^{k}|c_i|^2 B_i\|f\|^2, \  \ f\in L^2(\mathbb{R}).
	\end{align}
	By \eqref{st1} and \eqref{st2}, we conclude that  $\left\{ \sum\limits_{i=1}^{k} c_i\phi_{m,n,l}^{(i)^{P,Q}} \right\}_{m,n \in \mathbb{Z} \atop l \in \{1,2,\dots,N\}}$ forms a frame for $L^2(\mathbb{R})$ with the required frame bounds. This completes the proof.
	\endproof
Before an illustration of Theorem \ref{sfn}, we recall the following result.
	
	\begin{thm}\cite[Theorem 11.4.2]{ole}\label{te}
		Let $a$, $b>0$ and $\phi$ be a non-zero function in  $L^2(\mathbb{R})$. Suppose that
		\begin{align*}
			& B:=\frac{1}{b}\sup_{x\in[0,a]}\sum_{k\in \mathbb{Z}}\Big|\sum_{n\in \mathbb{Z}}\phi(x-na)\overline{\phi(x-na-k/b)}\Big|<\infty,
			\intertext{and}
			&A:=\frac{1}{b}\inf_{x\in[0,a]}\bigg(\sum_{n\in \mathbb{Z}}|\phi(x-na)|^2- \sum_{k \neq 0}\Big|\sum_{n \in \mathbb{Z}}\phi(x-na)\overline{\phi(x-na-k/b)}\Big|\bigg)> 0.
		\end{align*}
		Then, $\{e^{2\pi imbx}\phi(x-na)\}_{m,n\in \mathbb{Z}}$ is a frame for $ L^2(\mathbb{R})$ with frame  bounds $A$, $B$.
	\end{thm}
	
	The following example illustrates Theorem \ref{sfn}.
	\begin{exa}\label{finiteexa}
			Let $N=1$, $k=2$, and let $P = 1$, $Q =0$. Define a function $\phi^{(1)} := \phi \in L^2(\mathbb{R})$ as follows:
		\begin{align*}
			\phi (x)= \begin{cases}
				\sqrt{2x}, \quad &\text{if}\ x \in ]0,1];\\
				\sqrt{4x+2}, \quad &\text{if}\ x \in ]1,2]; \\
				0, \quad &\text{elsewhere}.
			\end{cases}
		\end{align*}
 Let $q_0=-1$, $p_0 =\pi$.  Then, $|p_0q_0|<2\pi$. Choose $a = 1$ and $b = \frac{1}{2}$ in Theorem \ref{te}.   Consider the  function $x \rightarrow \phi(x-n) \overline{\phi(x-n-2k)}$ for $x \in ]0,1]$ and  $n$, $k \in \mathbb{Z}$. It should be noted that due to the compact support of $\phi$, it can be nonzero only if $n \in \{-1,0\}$ and $k =0$. Consider the following functions:
		\begin{align*}
			G_0(x) = \sum_{n  \in \mathbb{Z}} \big|\phi(x-n)\big|^2 \quad \text{and} \quad G_1(x)= \sum_{k \neq 0} \Big|\sum_{n  \in \mathbb{Z}} \phi(x-n) \overline{\phi(x-n-2k)}\Big|, \ x \in ]0,1].
		\end{align*}
Then	
\begin{align*}
&G_0(x) = \sum_{n  \in \mathbb{Z}} \big|\phi(x-n)\big|^2 = \big|\phi(x)\big|^2 + \big|\phi(x+1)\big|^2 = 6x+2,  \ x \in ]0,1],
\intertext{and}
&G_1(x)= \sum_{k \neq 0} \Big|\sum_{n  \in \mathbb{Z}} \phi(x-n) \overline{\phi(x-n-2k)}\Big|=0, \ x \in ]0,1].
\end{align*}
This gives
		\begin{align*}
			A_1 = \frac{1}{b} \inf_{x\in]0,1]} \Big[ G_0(x) - G_1(x)\Big] = 4 \quad \text{and} \quad
			B_1 = \frac{1}{b} \sup_{x\in]0,1]} \Big[ G_0(x) + G_1(x)\Big] = 16.
		\end{align*}
		Therefore, by Theorem \ref{te}, $\big\{e^ {2 \pi i x \frac{m}{2}} \phi(x- n)\big\}_{m,n \in \mathbb{Z}}$ is a frame for $L^2(\mathbb{R})$ with frame bounds $4$ and $16$. That is
		\begin{align}\label{eqe1}
			4\|f\|^2 \leq \sum_{m,n \in \mathbb{Z}} \Big|\Big \langle f, e^{\pi i \langle m , \cdot \rangle}\phi(\cdot - n) \Big \rangle\Big|^2 \leq 16\|f\|^2, \ f\in L^2(\mathbb{R}).
		\end{align}
	For every $f \in L^2(\mathbb{R})$, we have
	\begin{align}\label{eqe2}
		\sum_{m, n\in \mathbb{Z}} \Big| \Big \langle f, \phi_{m,n}^{1,0}\Big \rangle \Big|^2 = \sum_{m, n\in \mathbb{Z}} \Big| \Big \langle f, e^{\frac{-\pi}{2}m n} e^{\pi i \langle m , \cdot \rangle}\phi(\cdot - n)\Big \rangle \Big|^2
		= \sum_{m,n \in \mathbb{Z}} \Big|\Big \langle f, e^{\pi i \langle m , \cdot \rangle}\phi(\cdot - n) \Big \rangle\Big|^2.
	\end{align}
	Using \eqref{eqe1} in \eqref{eqe2}, we conclude that  $\big\{\phi_{m,n}^{1,0}\big\}_{m, n\in\mathbb{Z}} = \big\{e^{\frac{-\pi}{2}m n} e^{\pi i \langle m , \cdot \rangle}\phi(\cdot - n)\big\}_{m, n\in\mathbb{Z}}$ is a frame for $L^2(\mathbb{R})$ with frame bounds $A_1=4$ and $B_1=16$.
		
Now consider the function $\phi^{(2)} := \psi \in L^2(\mathbb{R})$ which is given by:
		\begin{align*}
			\psi (x)= \begin{cases}
				2x, \quad &\text{if}\ x \in ]0,\frac{1}{2}];\\
				4(1-x), \quad &\text{if}\ x \in ]0,1];\\
				0, \quad &\text{elsewhere}.
			\end{cases}
		\end{align*}
	Let $q_0=\frac{1}{2}$, $p_0 =2\pi$.  Then, $|p_0q_0|<2\pi$. Choose $a = \frac{1}{2}$ and $b = 1$. Then , by Theorem \ref{te}, $\big\{e^ {2 \pi i x m} \psi\left(x-\frac{n}{2}\right)\big\}_{m,n \in \mathbb{Z}}$ is a frame for $L^2(\mathbb{R})$ with frame bounds $1$ and $4$. Hence, $\big\{\psi_{m,n}^{1,0}\big\}_{m, n\in\mathbb{Z}} =  \big\{e^{\frac{\pi}{2}m n} e^{2\pi i \langle m , \cdot \rangle}\psi\big(\cdot + \frac{n}{2}\big)\big\}_{m, n\in\mathbb{Z}}$ is a frame for $L^2(\mathbb{R})$ with frame bounds $A_2= 1$ and $B_2= 4$.
		
		Choose $c_1=-\frac{1}{2000}, c_2=\frac{1}{20}$ and $j=1$. Then, we compute inequality  given in \eqref{sfc}:
\begin{align*}
|c_1|A_1+ \sum_{i\neq 1}^{2}\bigg|\frac{c_i^2}{c_1}\bigg|A_1  = |c_1|A_1+\bigg|\frac{c_2^2}{c_1}\bigg|A_2 &= \frac{2501}{500}
 > \frac{4}{5} = 2|c_2| \sqrt{B_1B_2} = 2 \sqrt{B_j} \Big(\sum_{ i \neq 1}^{2}|c_i|\sqrt{B_i}\Big).
\end{align*}
Hence, by Theorem \ref{sfn}, the sum $\left\{c_1\phi_{m,n}^{1,0} + c_2\psi_{m,n}^{1,0}  \right\} _{m,n \in \mathbb{Z}}$ is a frame for $L^2(\mathbb{R})$ with frame bounds $\Big(\sum\limits_{i=1}^{2}|c_i|^2 A_i-2\sum\limits_{ i\neq 1}^{2}|c_1c_i|\sqrt{B_1 B_i}\Big) = 2101 \times 10^{-6}$ and $2\sum\limits_{i=1}^{2}|c_i|^2B_i=20008 \times 10^{-6}$.
	\end{exa}
	
 The authors of  \cite{OSCT} proved  that if $\{f_k\}_{k \in \mathcal{I}}$ is a frame for a separable Hilbert space $\mathcal{H}$  with frame operator $S$ and $\{g_k\}_{k \in \mathcal{I}}$ is an alternate dual frame then $\{S^a f_k +S^b g_k\}_{k \in \mathcal{I}}$ is a frame for $\mathcal{H}$ for all real numbers $a$, $b$.  They proved this result by using invertibility of an operator which is obtained by taking  sum of composition of   analysis operators of frames and bounded linear operators on the underlying  Hilbert space, see Corollary 3.3 of \cite{OSCT} for technical details. We show that sum of a frame and its dual frame  from the Weyl-Heisenberg group, with explicit frame bounds, constitutes a frame  by direct computation of both the lower and upper frame condition.
 \begin{thm}\label{dual}
		Let $\big\{\phi_{m,n,l}^{P, Q} \big\}_{m, n \in \mathbb{Z} \atop l \in \{1,2,\dots, N\} }$  be a frame for the space $L^2(\mathbb{R})$ with frame bounds $A_1$, $B_1$ and let $\big\{\psi_{m,n,l}^{P, Q} \big\}_{m, n \in \mathbb{Z} \atop l \in \{1,2,\dots, N\} }$ be a dual frame of $\big\{\phi_{m,n,l}^{P, Q} \big\}_{m, n \in \mathbb{Z} \atop l \in \{1,2,\dots, N\} }$ in $L^2(\mathbb{R})$ with frame bounds $A_2$, $B_2$. Then, $\big\{\phi_{m,n,l}^{P, Q}+ \psi_{m,n,l}^{P, Q} \big\}_{m, n \in \mathbb{Z} \atop  l \in \{1,2,\dots, N\} }$ is a frame for $L^2(\mathbb{R})$ with frame bounds $\big(A_1 + A_2 +2\big)$ and $\big(B_1 + B_2 +2\big)$.
	\end{thm}
	\proof
	For every $f \in L^2(\mathbb{R})$, we compute
	\begin{align}\label{dualeq2}
&\sum_{l = 1}^{N}\sum_{m,n\in\mathbb{Z}} \Big|\big \langle f, \phi_{m,n,l}^{P, Q}+ \psi_{m,n,l}^{P, Q} \big \rangle \Big|^2 \notag\\
 &= 	\sum_{l = 1}^{N}\sum_{m,n\in\mathbb{Z}} \big \langle f, \phi_{m,n,l}^{P, Q}+ \psi_{m,n,l}^{P, Q} \big \rangle \big \langle  \phi_{m,n,l}^{P, Q}+ \psi_{m,n,l}^{P, Q}, f \big \rangle \notag\\
		&= \sum_{l = 1}^{N}\sum_{m,n\in\mathbb{Z}} \Bigg(\Big(\langle f, \phi_{m,n,l}^{P, Q} \big \rangle + \langle f, \psi_{m,n,l}^{P, Q} \big \rangle\Big) \Big(\langle  \phi_{m,n,l}^{P, Q}, f \big \rangle + \langle  \psi_{m,n,l}^{P, Q}, f \big \rangle\Big)\Bigg) \notag\\
		&= \sum_{l = 1}^{N}\sum_{m,n\in\mathbb{Z}} \Bigg(\langle f, \phi_{m,n,l}^{P, Q} \big \rangle \langle  \phi_{m,n,l}^{P, Q}, f \big \rangle + \langle f, \phi_{m,n,l}^{P, Q} \big \rangle \langle  \psi_{m,n,l}^{P, Q}, f \big \rangle
 + \langle f, \psi_{m,n,l}^{P, Q} \big \rangle \langle  \phi_{m,n,l}^{P, Q}, f \big \rangle + \langle f, \psi_{m,n,l}^{P, Q} \big \rangle \langle  \psi_{m,n,l}^{P, Q}, f \big \rangle\Bigg) \notag\\
		&= \sum_{l = 1}^{N}\sum_{m,n\in\mathbb{Z}} \Big|\langle f, \phi_{m,n,l}^{P, Q} \big \rangle \Big|^2 + \Big \langle \sum_{l = 1}^{N}\sum_{m,n\in\mathbb{Z}}\big\langle f, \phi_{m,n,l}^{P, Q}\big \rangle \psi_{m,n,l}^{P, Q} , f \Big \rangle
		+ \Big \langle \sum_{l = 1}^{N}\sum_{m,n\in\mathbb{Z}}\big\langle f , \psi_{m,n,l}^{P, Q}\big \rangle \phi_{m,n,l}^{P, Q} , f \Big \rangle \notag\\
& \quad  + \sum_{l = 1}^{N}\sum_{m,n\in\mathbb{Z}} \Big|\langle f, \psi_{m,n,l}^{P, Q} \big \rangle \Big|^2 \notag\\
&= \sum_{l = 1}^{N}\sum_{m,n\in\mathbb{Z}} \Big|\langle f, \phi_{m,n,l}^{P, Q} \big \rangle \Big|^2 + \big\langle f , f \big \rangle + \big \langle f , f \big \rangle + \sum_{l = 1}^{N}\sum_{m,n\in\mathbb{Z}} \Big|\langle f, \psi_{m,n,l}^{P, Q} \big \rangle \Big|^2 \quad \Big(\text{by using} \ \ref{dualw} \Big) \notag \\
		& \leq B_1 \|f\|^2 + 2\|f\|^2 + B_2 \|f\|^2   \notag\\
		& = \Big(B_1 + B_2 +2\Big)\|f\|^2.
	\end{align}
	Similarly
	\begin{align}\label{dualeq3}
		\sum_{l = 1}^{N}\sum_{m,n\in\mathbb{Z}} \Big|\big \langle f, \phi_{m,n,l}^{P, Q}+ \psi_{m,n,l}^{P, Q} \big \rangle \Big|^2 &= \sum_{l = 1}^{N}\sum_{m,n\in\mathbb{Z}} \Big|\langle f, \phi_{m,n,l}^{P, Q} \big \rangle \Big|^2 + \big\langle f , f \big \rangle + \big \langle f , f \big \rangle + \sum_{l = 1}^{N}\sum_{m,n\in\mathbb{Z}} \Big|\langle f, \psi_{m,n,l}^{P, Q} \big \rangle \Big|^2 \notag \\
		& \geq A_1 \|f\|^2 + 2\|f\|^2 + A_2 \|f\|^2 \notag\\
		& = \Big(A_1 + A_2 +2\Big)\|f\|^2, \ f\in L^2(\mathbb{R}).
	\end{align}
The result now follows  from inequalities  \eqref{dualeq2} and \eqref{dualeq3}.
	\endproof
	
The following example is an illustration of Theorem \ref{dual}.
	\begin{exa}\label{dualexa}
		Let $N=2$ and $P=1$, $Q=0$. Define a function $\varphi_1 \in L^2(\mathbb{R})$ as
		\begin{align*}
			\varphi_1 (x)= \begin{cases}
				x-1, \quad &\text{if}\ x \in ]0,1];\\
				1-x, \quad &\text{if}\ x \in ]1,2];\\
				0, \quad &\text{elsewhere}.
			\end{cases}
		\end{align*}
	Let $q_0^{(1)}= -1$, $p_0^{(1)}=\pi$. Then, $|p_0^{(1)}q_0^{(1)}| <2\pi$. Choose  $a = 1$ and $b = \frac{1}{2}$ in Theorem \ref{te}.  Using the same method employed in Example \ref{bddexa}, it can be readily verified that $\big\{e^ {\pi i x m} \varphi_1(x- n)\big\}_{m,n \in \mathbb{Z}}$ is a tight frame for $L^2(\mathbb{R})$ with frame bound $2$. Therefore,
		\begin{align}\label{eqde1}
			\sum_{m,n \in \mathbb{Z}} \Big|\Big \langle f, e^{\pi i \langle m , \cdot \rangle}\varphi_1(\cdot - n) \Big \rangle\Big|^2 = 2\|f\|^2, \ f\in L^2(\mathbb{R}).
		\end{align}
		For $f \in L^2(\mathbb{R})$, using \eqref{eqde1}, we get
		\begin{align*}
			\sum_{m, n\in \mathbb{Z}} \Big| \Big \langle f, \varphi_{m,n,1}^{1,0}\Big \rangle \Big|^2 = \sum_{m, n\in \mathbb{Z}} \Big| \Big \langle f, e^{\frac{-\pi}{2}m n} e^{\pi i \langle m , \cdot \rangle}\varphi_1(\cdot - n)\Big \rangle \Big|^2
			= \sum_{m,n \in \mathbb{Z}} \Big|\Big \langle f, e^{\pi i \langle m , \cdot \rangle}\varphi_1(\cdot - n) \Big \rangle\Big|^2 =  2\|f\|^2.
		\end{align*}
	Hence, $\big\{\varphi_{m,n,1}^{1,0}\big\}_{m, n\in\mathbb{Z}} = \big\{e^{\frac{-\pi}{2}m n} e^{\pi i \langle m , \cdot \rangle}\varphi_1(\cdot - n)\big\}_{m, n\in\mathbb{Z}}$ forms a $2$-tight frame for $L^2(\mathbb{R})$. That is,
		\begin{align}\label{varphi1}
			\sum_{m,n\in\mathbb{Z}} \Big|\big \langle f, \varphi_{m,n,1}^{1,0} \big \rangle \Big|^2 = 2\|f\|^2,\ \text{for}\ f \in L^2(\mathbb{R}).
		\end{align}
		
		Similary, for the following function  $\varphi_2 \in L^2(\mathbb{R})$ as
		\begin{align*}
			\varphi_2 (x)= \begin{cases}
				\frac{x}{2}, \quad &\text{if}\ x \in ]0,1];\\
				\frac{1}{2}\sqrt{3-x}, \quad &\text{if}\ x \in ]1,2];\\
				0, \quad &\text{elsewhere},
			\end{cases}
		\end{align*}
		Let $q_0^{(2)}= -1$, $p_0^{(2)}=\pi$. Then, $|p_0^{(2)}q_0^{(2)}| <2\pi$. Choose $a = 1$ and $b = \frac{1}{2}$ in Theorem \ref{te}. Then, it is easy to observe that $\big\{e^ {\pi i x m} \varphi_2(x- n)\big\}_{m,n \in \mathbb{Z}}$ is a Parseval frame for $L^2(\mathbb{R})$. Hence, $\big\{\varphi_{m,n,2}^{1,0}\big\}_{m, n\in\mathbb{Z}} = \big\{e^{\frac{-\pi}{2}m n} e^{\pi i \langle m , \cdot \rangle}\varphi_2(\cdot - n)\big\}_{m, n\in\mathbb{Z}}$ is a Parseval frame for $L^2(\mathbb{R})$. That is,
		\begin{align}\label{varphi2}
			\sum_{m,n\in\mathbb{Z}} \Big|\big \langle f, \varphi_{m,n,2}^{1,0} \big \rangle \Big|^2 = \|f\|^2,\ \text{for}\ f \in L^2(\mathbb{R}).
		\end{align}
		Thus, for every $f \in L^2(\mathbb{R})$, we have
		\begin{align}\label{varphi}
			\sum_{l = 1}^{N} \sum_{m,n\in\mathbb{Z}} \Big|\big \langle f, \varphi_{m,n,l}^{1,0} \big \rangle \Big|^2= \sum_{m,n\in\mathbb{Z}} \Big|\big \langle f, \varphi_{m,n,1}^{1,0} \big \rangle \Big|^2 + \sum_{m,n\in\mathbb{Z}} \Big|\big \langle f, \varphi_{m,n,2}^{1,0} \big \rangle \Big|^2 = 3\|f\|^2.
		\end{align}
		Hence, the sequence $\big\{\varphi_{m,n,l}^{1,0}\big\}_{m, n\in\mathbb{Z} \atop l\in\{1,2\}}$ is a $3$-tight frame for $L^2(\mathbb{R})$. Therefore, by Theorem \ref{tightthm}, every $f \in L^2(\mathbb{R})$ can be written as 	
		\begin{align}\label{dualexa1}
			f = \frac{1}{3}\sum_{l = 1}^{N} \sum_{m,n\in\mathbb{Z}} \big \langle  f, \varphi_{m,n,l}^{1,0} \big \rangle \varphi_{m,n,l}^{1,0}.
		\end{align}
		
		Define $\big\{\phi_{m,n,l}^{1,0} \big\}_{m, n \in \mathcal{I} \atop l \in \{1,2\} }= \big\{\varphi_{m,n,1}^{1,0}\big\}_{m,n<0} \cup \big\{\frac{1}{3}\varphi_{m,n,l}^{1,0} \big\}_{m, n \in \mathbb{Z} \atop l \in \{1,2\} }$, where $\mathcal{I}$ is any countable set. Then, for $f \in L^2(\mathbb{R})$, we obtain
		\begin{align}\label{phi1}
			\sum_{l = 1}^{2}\sum_{m,n\in \mathcal{I}} \Big|\big \langle f, \phi_{m,n,l}^{1,0} \big \rangle \Big|^2 & = \sum_{m,n<0} \Big|\big \langle f, \varphi_{m,n,1}^{1,0} \big \rangle \Big|^2+ \frac{1}{9}\sum_{l = 1}^{2}\sum_{m,n\in\mathbb{Z}} \Big|\big \langle f, \varphi_{m,n,l}^{1,0} \big \rangle \Big|^2 \notag\\
			& \leq \sum_{m,n \in \mathbb{Z}} \Big|\big \langle f, \varphi_{m,n,1}^{1,0} \big \rangle \Big|^2+ \frac{1}{9}\sum_{l = 1}^{2}\sum_{m,n\in\mathbb{Z}} \Big|\big \langle f, \varphi_{m,n,l}^{1,0} \big \rangle \Big|^2.
		\end{align}
		Using \eqref{varphi1} and \eqref{varphi} in \eqref{phi1}, we get
		\begin{align}\label{phi2}
			\sum_{l = 1}^{2}\sum_{m,n\in \mathcal{I}} \Big|\big \langle f, \phi_{m,n,l}^{1,0} \big \rangle \Big|^2 \leq 2\|f\|^2 +\frac{3}{9}\|f\|^2 \ = \frac{7}{3}\|f\|^2, \ f \in L^2(\mathbb{R}).
		\end{align}
		Similarly, using \eqref{varphi1}, we have
		\begin{align}\label{phi3}
			\sum_{l = 1}^{2}\sum_{m,n\in \mathcal{I}} \Big|\big \langle f, \phi_{m,n,l}^{1,0} \big \rangle \Big|^2 & = \sum_{m,n<0} \Big|\big \langle f, \varphi_{m,n,1}^{1,0} \big \rangle \Big|^2+ \frac{1}{9}\sum_{l = 1}^{2}\sum_{m,n\in\mathbb{Z}} \Big|\big \langle f, \varphi_{m,n,l}^{1,0} \big \rangle \Big|^2 \notag\\
			&\geq  \frac{1}{9}\sum_{l = 1}^{2}\sum_{m,n\in\mathbb{Z}} \Big|\big \langle f, \varphi_{m,n,l}^{1,0} \big \rangle \Big|^2 = \frac{1}{3}\|f\|^2, \ f \in L^2(\mathbb{R}).	
		\end{align}
		Hence, from \eqref{phi2} and \eqref{phi3}, the sequence $\big\{\phi_{m,n,l}^{1,0} \big\}_{m, n \in \mathcal{I} \atop l \in \{1,2\}}$ is a frame for $ L^2(\mathbb{R})$ with frame bounds $A_1= \frac{1}{3}$ and $B_1 = \frac{7}{3}$.
		
		Now, define the sequence $\big\{\psi_{m,n,l}^{1,0} \big\}_{m, n \in \mathcal{I} \atop l \in \{1,2\} }= \big\{\frac{1}{6}\varphi_{m,n,1}^{1,0}\big\}_{m,n<0} \cup \big\{\frac{1}{2}\varphi_{m,n,1}^{1,0}\big\}_{m,n<0} \cup \big\{\varphi_{m,n,1}^{1,0}\big\}_{m,n\geq 0} \cup  \big\{\varphi_{m,n,2}^{1,0} \big\}_{m, n \in \mathbb{Z}}$. Then,  for every $f \in L^2(\mathbb{R})$, we have
		\begin{align*}
			\sum_{l = 1}^{2}\sum_{m,n\in \mathcal{I}} \big \langle  f, \phi_{m,n,l}^{1,0} \big \rangle \psi_{m,n,l}^{1,0} &= \frac{1}{6}\sum_{m,n<0} \big \langle  f, \varphi_{m,n,1}^{1,0} \big \rangle \varphi_{m,n,1}^{1,0} + \frac{1}{6} \sum_{m,n<0} \big \langle  f, \varphi_{m,n,1}^{1,0}\big \rangle \varphi_{m,n,1}^{1,0} \\
			& + \frac{1}{3} \sum_{m,n\geq 0} \big \langle  f, \varphi_{m,n,1}^{1,0} \big \rangle \varphi_{m,n,1}^{1,0} + \frac{1}{3} \sum_{m,n \in \mathbb{Z}} \big \langle  f, \varphi_{m,n,2}^{1,0} \big \rangle \varphi_{m,n,2}^{1,0}\\
			& = \frac{1}{3} \sum_{l = 1}^{2}\sum_{m,n\in\mathbb{Z}} \big \langle f,\varphi_{m,n,l}^{1,0}\big \rangle \varphi_{m,n,l}^{1,0} = f  \quad  \text{\big(using \eqref{dualexa1}\big)}.
		\end{align*}
		Therefore, $\big\{\psi_{m,n,l}^{1,0} \big\}_{m, n \in \mathcal{I} \atop l \in \{1,2\} }$ is a dual frame of $\big\{\phi_{m,n,l}^{1,0} \big\}_{m, n \in \mathcal{I} \atop l \in \{1,2\}}$. Now, we will evaluate frame bounds of $\big\{\psi_{m,n,l}^{1,0} \big\}_{m, n \in \mathcal{I} \atop l \in \{1,2\} }$. Using \eqref{varphi1} and \eqref{varphi2}, for each $f \in L^2(\mathbb{R})$,
		\begin{align}\label{psi1}
			&\sum_{l = 1}^{2}\sum_{m,n\in \mathcal{I}} \Big|\big \langle f, \psi_{m,n,l}^{1,0} \big \rangle \Big|^2  \notag\\
& = \frac{1}{36} \sum_{m,n <0} \Big|\big \langle f, \varphi_{m,n,1}^{1,0} \big \rangle \Big|^2 + \frac{1}{4} \sum_{m,n <0} \Big|\big \langle f, \varphi_{m,n,1}^{1,0} \big \rangle \Big|^2
			 + \sum_{m,n \geq 0} \Big|\big \langle f, \varphi_{m,n,1}^{1,0} \big \rangle \Big|^2 + \sum_{m,n \in \mathbb{Z}} \Big|\big \langle f, \varphi_{m,n,2}^{1,0} \big \rangle \Big|^2 \notag\\
			& \leq  \sum_{m,n <0} \Big|\big \langle f, \varphi_{m,n,1}^{1,0} \big \rangle \Big|^2 +  \sum_{m,n <0} \Big|\big \langle f, \varphi_{m,n,1}^{1,0} \big \rangle \Big|^2   + \sum_{m,n \geq 0} \Big|\big \langle f, \varphi_{m,n,1}^{1,0} \big \rangle \Big|^2 + \sum_{m,n \in \mathbb{Z}} \Big|\big \langle f, \varphi_{m,n,2}^{1,0} \big \rangle \Big|^2\notag \\
			& \leq 2 \sum_{m,n \in \mathbb{Z}} \Big|\big \langle f, \varphi_{m,n,1}^{1,0} \big \rangle \Big|^2 + \sum_{m,n \in \mathbb{Z}} \Big|\big \langle f, \varphi_{m,n,2}^{1,0} \big \rangle \Big|^2\notag \\
			& = 4\|f\|^2 + \|f\|^2 = 5 \|f\|^2.
		\end{align}
		Similarly, using \eqref{varphi2},
		\begin{align}\label{psi2}
			\sum_{l = 1}^{2}\sum_{m,n\in I} \Big|\big \langle f, \psi_{m,n,l}^{1,0} \big \rangle \Big|^2 & \ge \sum_{m,n \in I} \Big|\big \langle f, \varphi_{m,n,2}^{1,0} \big \rangle \Big|^2 = \|f\|^2, \ f \in L^2(\mathbb{R}).
		\end{align}
		Hence, $\big\{\psi_{m,n,l}^{1,0} \big\}_{m, n \in \mathcal{I} \atop l \in \{1,2\} }$ is a frame with frame bounds $A_2 = 1$ and $B_2 =5$. Thus, by Theorem \ref{dual}, $\big\{\phi_{m,n,l}^{1,0}+ \psi_{m,n,l}^{1,0}\big\}_{m, n \in \mathcal{I} \atop  l \in \{1,2\} }$ is a frame with frame bounds $\frac{10}{3}$ and $\frac{28}{3}$.
\end{exa}

It is well known that the image of a frame for a space  under a bounded linear operator may not be a frame for the underlying space. In this direction, we give sufficient conditions for  sum of images of frames under bounded linear operators in terms of lower bounds of Hilbert-adjoint of bounded linear  operators and frame bounds of given frames.
\begin{thm}\label{operator}
		Let $\big\{\phi_{m,n,l}^{P, Q} \big\}_{m, n \in \mathbb{Z} \atop l \in \{1,2,\dots, N\} }$ and $\big\{\psi_{m,n,l}^{P, Q}\big\}_{m, n \in \mathbb{Z} \atop l \in \{1,2,\dots, N\} }$ be  frames for $L^2(\mathbb{R})$ with frame bounds $A_1$, $B_1$ and $A_2$, $B_2$, respectively. Let  $\Theta_1$, $\Theta_2$ be  bounded linear operators acting on  $L^2(\mathbb{R})$ such that
\begin{enumerate}[$(i)$]
  \item $\Theta_1^*$, $\Theta_2^*$ are bounded below by  constants $m_1$ and $m_2$.
  \item $A_1m_1^2 + A_2 m_2^2 > 2\sqrt{B_1 B_2}\|\Theta_1\|\|\Theta_2\|$.
\end{enumerate}
 Then, $\Big\{\Theta_1 \phi_{m,n,l}^{P, Q}+ \Theta_2 \psi_{m,n,l}^{P, Q}\Big\}_{m, n \in \mathbb{Z} \atop  l \in \{1,2,\dots, N\} }$ is a frame for $L^2(\mathbb{R})$ with frame bounds
		\begin{align*}
\Big(A_1m_1^2 + A_2 m_2^2 - 2\sqrt{B_1 B_2}\|\Theta_1\|\|\Theta_2\|\Big) \quad  \text{and} \quad  \Big(\sqrt{B_1}\|\Theta_1\|+ \sqrt{B_2}\|\Theta_2\|\Big)^2.
		\end{align*}
	\end{thm}	
	\proof	
For every $f \in L^2(\mathbb{R})$, we compute
	\begin{align}
&\sum_{l = 1}^{N}\sum_{m,n\in\mathbb{Z}} \Big|\big \langle f, \Theta_1\phi_{m,n,l}^{P, Q}+ \Theta_2 \psi_{m,n,l}^{P, Q} \big \rangle \Big|^2 \notag\\
&= 	\sum_{l = 1}^{N}\sum_{m,n\in\mathbb{Z}} \big \langle f, \Theta_1\phi_{m,n,l}^{P, Q}+ \Theta_2 \psi_{m,n,l}^{P, Q} \big \rangle \big \langle  \Theta_1\phi_{m,n,l}^{P, Q}+ \Theta_2 \psi_{m,n,l}^{P, Q}, f \big \rangle \notag\\
		&= \sum_{l = 1}^{N}\sum_{m,n\in\mathbb{Z}} \Big(\big\langle f, \Theta_1\phi_{m,n,l}^{P, Q} \big \rangle + \big\langle f, \Theta_2 \psi_{m,n,l}^{P, Q} \big \rangle\Big) \Big(\big\langle  \Theta_1\phi_{m,n,l}^{P, Q}, f \big \rangle +  \big\langle  \Theta_2 \psi_{m,n,l}^{P, Q}, f \big \rangle\Big)\notag\\
		&= \sum_{l = 1}^{N}\sum_{m,n\in\mathbb{Z}} \Big| \big\langle f, \Theta_1\phi_{m,n,l}^{P, Q} \big \rangle \Big|^2 + \sum_{l = 1}^{N}\sum_{m,n\in\mathbb{Z}} \big \langle f, \Theta_1\phi_{m,n,l}^{P, Q} \big \rangle  \big\langle  \Theta_2 \psi_{m,n,l}^{P, Q}, f \big \rangle \notag\\
		&  + \sum_{l = 1}^{N}\sum_{m,n\in\mathbb{Z}}\big \langle f, \Theta_2 \psi_{m,n,l}^{P, Q} \big \rangle \big\langle  \Theta_1\phi_{m,n,l}^{P, Q}, f \big \rangle + \sum_{l = 1}^{N}\sum_{m,n\in\mathbb{Z}} \Big|\big\langle f, \Theta_2 \psi_{m,n,l}^{P, Q}\big \rangle \Big|^2 \notag\\
		&= \sum_{l = 1}^{N}\sum_{m,n\in\mathbb{Z}} \Big|\big\langle \Theta_1^*f, \phi_{m,n,l}^{P, Q} \big \rangle \Big|^2 +  \sum_{l = 1}^{N}\sum_{m,n\in\mathbb{Z}}\big \langle f, \Theta_1\phi_{m,n,l}^{P, Q} \big \rangle \langle   \Theta_2\psi_{m,n,l}^{P, Q}, f \big \rangle \notag\\
		& +\quad \overline{ \sum_{l = 1}^{N}\sum_{m,n\in\mathbb{Z}} \big\langle f, \Theta_1\phi_{m,n,l}^{P, Q} \big \rangle \big \langle   \Theta_2\psi_{m,n,l}^{P, Q}, f \big \rangle} + \sum_{l = 1}^{N}\sum_{m,n\in\mathbb{Z}} \Big|\big \langle \Theta_2^*  f, \psi_{m,n,l}^{P, Q}\big \rangle \Big|^2\notag\\
		&=  \sum_{l = 1}^{N}\sum_{m,n\in\mathbb{Z}} \Big|\big\langle \Theta_1^*f, \phi_{m,n,l}^{P, Q} \big \rangle \Big|^2 + \sum_{l = 1}^{N}\sum_{m,n\in\mathbb{Z}} \Big|\big\langle \Theta_2^*  f, \psi_{m,n,l}^{P, Q}\big \rangle \Big|^2    + 2 \mathrm{Re}  \sum_{l = 1}^{N}\sum_{m,n\in\mathbb{Z}} \big\langle f, \Theta_1\phi_{m,n,l}^{P, Q} \big \rangle \big\langle   \Theta_2\psi_{m,n,l}^{P, Q}, f \big \rangle\notag\\
		& \leq \sum_{l = 1}^{N}\sum_{m,n\in\mathbb{Z}} \Big|\big\langle \Theta_1^*f, \phi_{m,n,l}^{P, Q} \big \rangle \Big|^2 + \sum_{l = 1}^{N}\sum_{m,n\in\mathbb{Z}} \Big|\big\langle \Theta_2^*  f, \psi_{m,n,l}^{P, Q}\big \rangle \Big|^2   +2 \Big| \sum_{l = 1}^{N}\sum_{m,n\in\mathbb{Z}} \big\langle f, \Theta_1\phi_{m,n,l}^{P, Q} \big \rangle \big\langle   \Theta_2\psi_{m,n,l}^{P, Q}, f \big \rangle\Big|\notag\\
		& \leq \sum_{l = 1}^{N}\sum_{m,n\in\mathbb{Z}} \Big|\big\langle \Theta_1^*f, \phi_{m,n,l}^{P, Q} \big \rangle \Big|^2 + \sum_{l = 1}^{N}\sum_{m,n\in\mathbb{Z}} \Big|\big\langle \Theta_2^*  f, \psi_{m,n,l}^{P, Q}\big \rangle \Big|^2   +2 \sum_{l = 1}^{N}\sum_{m,n\in\mathbb{Z}} \Big|\big\langle \Theta_1^*f, \phi_{m,n,l}^{P, Q} \big \rangle \Big| \
\Big|\big \langle   \psi_{m,n,l}^{P, Q}, \Theta_2^* f \big \rangle \Big| \notag\\
&\leq \sum_{l = 1}^{N}\sum_{m,n\in\mathbb{Z}} \Big|\big\langle \Theta_1^*f, \phi_{m,n,l}^{P, Q} \big \rangle \Big|^2 + \sum_{l = 1}^{N}\sum_{m,n\in\mathbb{Z}} \Big|\big\langle \Theta_2^*  f, \psi_{m,n,l}^{P, Q}\big \rangle \Big|^2 \notag\\
& + 2\bigg(\sum_{l = 1}^{N}\sum_{m,n\in\mathbb{Z}} \Big|\big\langle \Theta_1^*f, \phi_{m,n,l}^{P, Q} \big \rangle \Big|^2\bigg)^{\frac{1}{2}} \Bigg(\sum_{l = 1}^{N}\sum_{m,n\in\mathbb{Z}} \Big|\big\langle \Theta_2^*  f, \psi_{m,n,l}^{P, Q}\big \rangle \Big|^2\Bigg)^{\frac{1}{2}}\notag\\
		& \leq B_1\|\Theta_1^*f\|^2  + B_2 \|\Theta_2^* f\|^2  + 2\Big(\sqrt{B_1}\|\Theta_1^*f\|\Big) \Big(\sqrt{B_2}\|\Theta_2^* f\|\Big)\notag \\
		&\leq  B_1\|\Theta_1^*\|^2\|f\|^2  + B_2 \|\Theta_2^*\|^2\| f\|^2  + 2\sqrt{B_1 B_2} \|\Theta_1^*\| \|\Theta_2^*\| \| f\|^2 \notag \\
                  &=  B_1\|\Theta_1\|^2\|f\|^2  + B_2 \|\Theta_2\|^2\| f\|^2  + 2\sqrt{B_1 B_2} \|\Theta_1\| \|\Theta_2\| \| f\|^2 \notag \\
		&= \Big(\sqrt{B_1} \|\Theta_1\|+ \sqrt{B_2}  \|\Theta_2\|\Big)^2\|f\|^2, \ f \in L^2(\mathbb{R})\notag.
	\end{align}
This gives the required upper frame bound for the collection $\Big\{\Theta_1 \phi_{m,n,l}^{P, Q}+ \Theta_2 \psi_{m,n,l}^{P, Q}\Big\}_{m, n \in \mathbb{Z} \atop  l \in \{1,2,\dots, N\} }$.

For the lower frame bound, by employing the inequality $|x+y|^2 \geq \Big(|x|-|y|\Big)^2$ for any $x,y \in \mathbb{C}$, we compute
	\begin{align*}
&\sum_{l = 1}^{N}\sum_{m,n\in\mathbb{Z}} \Big|\big \langle f, \Theta_1\phi_{m,n,l}^{P, Q}+ \Theta_2 \psi_{m,n,l}^{P, Q} \big \rangle \Big|^2\\
 &= \sum_{l = 1}^{N}\sum_{m,n\in\mathbb{Z}} \Big|\big \langle f, \Theta_1\phi_{m,n,l}^{P, Q}\big \rangle + \big\langle f, \Theta_2 \psi_{m,n,l}^{P, Q} \big \rangle\Big|^2 \\
		& \geq \sum_{l = 1}^{N}\sum_{m,n\in\mathbb{Z}} \Big|\big \langle f, \Theta_1\phi_{m,n,l}^{P, Q} \big \rangle \Big|^2 + \sum_{l = 1}^{N}\sum_{m,n\in\mathbb{Z}} \Big|\big \langle f, \Theta_2 \psi_{m,n,l}^{P, Q} \big \rangle \Big|^2   -2 \sum_{l = 1}^{N}\sum_{m,n\in\mathbb{Z}} \Big|\big \langle f, \Theta_1\phi_{m,n,l}^{P, Q}\big \rangle \Big| \  \Big|\big \langle f, \Theta_2 \psi_{m,n,l}^{P, Q} \big \rangle \Big| \\
		& \geq \sum_{l = 1}^{N}\sum_{m,n\in\mathbb{Z}} \Big|\big \langle \Theta_1^*f, \phi_{m,n,l}^{P, Q} \big \rangle \Big|^2 + \sum_{l = 1}^{N}\sum_{m,n\in\mathbb{Z}} \Big|\big \langle \Theta_2^* f, \psi_{m,n,l}^{P, Q} \big \rangle \Big|^2  \notag \\
& - 2\bigg(\sum_{l = 1}^{N}\sum_{m,n\in\mathbb{Z}} \Big|\big\langle \Theta_1^* f, \phi_{m,n,l}^{P, Q} \big \rangle \Big|^2\bigg)^{\frac{1}{2}} \Bigg(\sum_{l = 1}^{N}\sum_{m,n\in\mathbb{Z}} \Big|\big\langle \Theta_2^*  f, \phi_{m,n,l}^{P, Q}\big \rangle \Big|^2\Bigg)^{\frac{1}{2}}\\
		&\geq A_1\|\Theta_1^* f\|^2+ A_2 \|\Theta_2^* f\|^2  -2 \Big(\sqrt{B_1}\|\Theta_1^*f\| \Big) \Big(\sqrt{B_2}\|\Theta_2^*f\|\Big)\notag\\
		& \geq A_1\|\Theta_1^* f\|^2+ A_2 \|\Theta_2^* f\|^2 - 2\sqrt{B_1 B_2}\|\Theta_1\| \|\Theta_2\|\|f\|^2 \\
& \geq A_1 m_1^2\|f\|^2 + A_2 m_2^2\|f\|^2 - 2\sqrt{B_1 B_2}\|\Theta_1\| \|\Theta_2\|\|f\|^2 \\
		& = \Big( A_1 m_1^2 + A_2 m_2^2 - 2\sqrt{B_1 B_2}\|\Theta_1\| \|\Theta_2\|\Big)\|f\|^2, \ f \in L^2(\mathbb{R}).
	\end{align*}
This completes the proof.
	\endproof
\begin{rem}
The authors in  \cite[Proposition 3.1]{OSCT} proved equivalent conditions for the sum of images of Bessel sequences under bounded linear operator to be a frame in terms of invertibility of an operator which is obtained by taking sum of composition of analysis operator of Bessel sequences and bounded linear operator on the underlying Hilbert space. Conditions in our result are different than that given in \cite[Proposition 3.1]{OSCT}.
\end{rem}
The following example illustrates Theorem \ref{operator}.
\begin{exa}\label{operexa1}
	Let $N=1$ and let $\big\{\phi_{m,n}^{1,0}\big\}_{m, n\in\mathbb{Z}}$ be a frame for $L^2(\mathbb{R})$ given in Example \ref{finiteexa} with frame bounds $A_1=4$ and $B_1=16$. Define a  function $\psi \in L^2(\mathbb{R})$ by
	\begin{align*}
		\psi (x)= \begin{cases}
			\sqrt{2}x, \quad &\text{if}\ x \in ]0,1];\\
			\sqrt{2} x -2\sqrt{2}, \quad &\text{if}\ x \in ]1,2];\\
			0, \quad &\text{elsewhere}.
		\end{cases}
	\end{align*}
Let $q_0= 1$, $p_0= \pi$, $a = 1$ and $b = \frac{1}{2}$. By  invoking Theorem \ref{te},   $\big\{e^ {2 \pi i x \frac{m}{2}} \psi(x- n)\big\}_{m,n \in \mathbb{Z}}$ is a $4$-tight frame for $L^2(\mathbb{R})$. Hence, $\big\{\psi_{m,n}^{1,0}\big\}_{m, n\in\mathbb{Z}}$ is a tight frame for $L^2(\mathbb{R})$ with frame bound $A_2=B_2 =4$.
	
	Define bounded linear operators $\Theta_1$, $\Theta_2 \colon L^2(\mathbb{R}) \rightarrow L^2(\mathbb{R})$ by
	\begin{align*}
		\Theta_1f(x) =\frac{1}{4} f(x) \quad \text{and} \quad
		\Theta_2 f(x)= f(x), \ f \in L^2(\mathbb{R}).
	\end{align*}
Then, $\Theta_1$ and $\Theta_2$ are bounded linear operators such that $\|\Theta_1\| = \frac{1}{4}$, $\|\Theta_2\| =1$. Also,  $\|\Theta_1^* f\| = \frac{1}{4}\| f\|$ and  $  \|\Theta_2^* f\| = 1 \|f\|$. So, we can take $m_1 = \frac{1}{4}$ and $m_2 =1$.
Therefore
	\begin{align*}
	A_1m_1^2 + A_2 m_2^2 = 4 \cdot \frac{1}{16} + 4 \cdot 1 = \frac{17}{4} > 4 = 2\sqrt{B_1 B_2}\|\Theta_1\|\|\Theta_2\|.
	\end{align*}
	Hence, by Theorem \ref{operator}, $\big\{\Theta_1 \phi_{m,n}^{P, Q}+ \Theta_2 \psi_{m,n}^{P, Q} \big\}_{m, n \in \mathbb{Z}}$ is a frame for $L^2(\mathbb{R})$ with frame bounds $\frac{1}{4}$ and $9$.
\end{exa}

	The following theorem establishes sufficient conditions on bounded sequences of scalars in $\mathbb{C}$ and provides a relationship between them and the bounds of frames $\big\{\phi_{m,n,l}^{P, Q} \big\}_{m, n \in \mathbb{Z} \atop l \in {1,2,\dots, N} }$ and $\big\{\psi_{m,n,l}^{P, Q} \big\}_{m, n \in \mathbb{Z} \atop l \in {1,2,\dots, N} }$ associated with the Weyl-Heisenberg group in $L^2(\mathbb{R})$. The theorem also proves that the sum of these two frames forms a frame.
	
	\begin{thm}\label{bddseq}
Let $\big\{\phi_{m,n,l}^{P, Q} \big\}_{m, n \in \mathbb{Z} \atop l \in \{1,2,\dots, N\} }$ and  $\big\{\psi_{m,n,l}^{P, Q}\big\}_{m, n \in \mathbb{Z} \atop l \in \{1,2,\dots, N\} }$ be frames for $L^2(\mathbb{R})$ with frame bounds $A_1$, $B_1$;  and $A_2$, $B_2$, respectively. Let $\big\{\alpha_{m,n,l} \big\}_{m, n \in \mathbb{Z} \atop l \in \{1,2,\dots, N\} }$ and  $\big\{\beta_{m,n,l} \big\}_{m, n \in \mathbb{Z} \atop l \in \{1,2,\dots, N\} }$ be bounded sequences of scalars such that
		\begin{align*}
			2 \sup |\alpha_{m,n,l}| \sup |\beta_{m,n,l}| \sqrt{B_1 B_2} < \inf |\alpha_{m,n,l}|^2 A_1 + \inf |\beta_{m,n,l}|^2 A_2.
		\end{align*}
	 Then, the perturbed sum $\Big\{\alpha_{m,n,l}\phi_{m,n,l}^{P, Q}+ \beta_{m,n,l}\psi_{m,n,l}^{P, Q}\Big\}_{m, n \in \mathbb{Z} \atop  l \in \{1,2,\dots, N\} }$ is a frame for $L^2(\mathbb{R})$ with  frame bounds
		\begin{align*}
&\Big(\inf |\alpha_{m,n,l}|^2 A_1 + \inf |\beta_{m,n,l}|^2 A_2 - 2 \sup |\alpha_{m,n,l}| \sup |\beta_{m,n,l}| \sqrt{B_1 B_2}\Big)
 \intertext{and}
& \Big(\sup |\alpha_{m,n,l}| \sqrt{B_1} + \sup |\beta_{m,n,l}| \sqrt{B_2}\Big)^2.
		\end{align*}
	\end{thm}
	\proof
For every  $f \in L^2(\mathbb{R})$, we compute
	\begin{align}\label{bddeq1}
		&\sum_{l = 1}^{N}\sum_{m,n\in\mathbb{Z}} \Big|\big \langle f, \alpha_{m,n,l}\phi_{m,n,l}^{P, Q}+ \beta_{m,n,l}\psi_{m,n,l}^{P, Q} \big \rangle \Big|^2 \notag\\
		& \geq \sum_{l = 1}^{N}\sum_{m,n\in\mathbb{Z}} \Big|\big \langle f, \alpha_{m,n,l}\phi_{m,n,l}^{P, Q} \big \rangle \Big|^2  + \sum_{l = 1}^{N}\sum_{m,n\in\mathbb{Z}} \Big|\big \langle f, \beta_{m,n,l}\psi_{m,n,l}^{P, Q} \big \rangle \Big|^2  \notag\\
		&  -2 \sum_{l = 1}^{N}\sum_{m,n\in\mathbb{Z}} \Big|\big \langle f, \alpha_{m,n,l}\phi_{m,n,l}^{P, Q} \big \rangle \Big| \Big|\big \langle f,  \beta_{m,n,l}\psi_{m,n,l}^{P, Q} \big \rangle \Big| \notag\\
		&\geq \sum_{l = 1}^{N}\sum_{m,n\in\mathbb{Z}}\big|\alpha_{m,n,l}\big|^2 \Big|\big \langle f, \phi_{m,n,l}^{P, Q} \big \rangle \Big|^2  +  \sum_{l = 1}^{N}\sum_{m,n\in\mathbb{Z}} \big|\beta_{m,n,l}\big|^2\Big|\big \langle f, \psi_{m,n,l}^{P, Q} \big \rangle \Big|^2 \notag\\
		& -2 \Bigg(\sum_{l = 1}^{N}\sum_{m,n\in\mathbb{Z}} \Big|\big \langle f, \alpha_{m,n,l}\phi_{m,n,l}^{P, Q} \big \rangle \Big|^2\Bigg)^{\frac{1}{2}} \Bigg(\sum_{l = 1}^{N}\sum_{m,n\in\mathbb{Z}} \Big|\big \langle f,  \beta_{m,n,l}\psi_{m,n,l}^{P, Q} \big \rangle \Big|^2 \Bigg)^{\frac{1}{2}}\notag\\
		& \geq \inf |\alpha_{m,n,l}|^2 \sum_{l = 1}^{N}\sum_{m,n\in\mathbb{Z}} \Big|\big \langle f, \phi_{m,n,l}^{P, Q} \big \rangle \Big|^2 + \inf |\beta_{m,n,l}|^2 \sum_{l = 1}^{N}\sum_{m,n\in\mathbb{Z}} \Big|\big \langle f, \phi_{m,n,l}^{P, Q} \big \rangle \Big|^2 \notag\\
		&  -2\sup |\alpha_{m,n,l}|\sup |\beta_{m,n,l}| \Bigg(\sum_{l = 1}^{N}\sum_{m,n\in\mathbb{Z}} \Big|\big \langle f, \phi_{m,n,l}^{P, Q} \big \rangle \Big|^2\Bigg)^{\frac{1}{2}} \Bigg(\sum_{l = 1}^{N}\sum_{m,n\in\mathbb{Z}} \Big|\big \langle f, \psi_{m,n,l}^{P, Q} \big \rangle \Big|^2 \Bigg)^{\frac{1}{2}}\notag\\
		&\geq \inf |\alpha_{m,n,l}|^2 A_1\|f\|^2 + \inf |\beta_{m,n,l}|^2 A_2\|f\|^2 - 2 \sup |\alpha_{m,n,l}| \sup |\beta_{m,n,l}| \sqrt{B_1 B_2}\|f\|^2 \notag\\
		&= \Big(\inf |\alpha_{m,n,l}|^2 A_1 + \inf |\beta_{m,n,l}|^2 A_2 - 2 \sup |\alpha_{m,n,l}| \sup |\beta_{m,n,l}| \sqrt{B_1 B_2}\Big) \|f\|^2.
	\end{align}
Now, for any $f \in L^2(\mathbb{R})$, we have
	\begin{align*}
		&\Bigg(\sum_{l = 1}^{N}\sum_{m,n\in\mathbb{Z}} \Big|\big \langle f, \alpha_{m,n,l}\phi_{m,n,l}^{P, Q}+ \beta_{m,n,l}\psi_{m,n,l}^{P, Q} \big \rangle \Big|^2 \Bigg)^{\frac{1}{2}}\\
		& \leq \Bigg(\sum_{l = 1}^{N}\sum_{m,n\in\mathbb{Z}} \Big|\big \langle f, \alpha_{m,n,l}\phi_{m,n,l}^{P, Q} \big \rangle \Big|^2 \Bigg)^{\frac{1}{2}} + \Bigg(\sum_{l = 1}^{N}\sum_{m,n\in\mathbb{Z}} \Big|\big \langle f,  \beta_{m,n,l}\psi_{m,n,l}^{P, Q} \big \rangle \Big|^2 \Bigg)^{\frac{1}{2}} \notag\\
		& \leq \sup \ |\alpha_{m,n,l}| \sqrt{B_1}\|f\| + \sup \ |\beta_{m,n,l}| \sqrt{B_2}\|f\|,
	\end{align*}
which entails
	\begin{align}\label{bddeq2}
		\sum_{l = 1}^{N}\sum_{m,n\in\mathbb{Z}} \Big|\big \langle f, \alpha_{m,n,l}\phi_{m,n,l}^{P, Q}+ \beta_{m,n,l}\psi_{m,n,l}^{P, Q} \big \rangle \Big|^2 \leq \Big(\sup \ |\alpha_{m,n,l}| \sqrt{B_1} + \sup \ |\beta_{m,n,l}| \sqrt{B_2}\Big)^2\|f\|^2
	\end{align}
	for all $f \in L^2(\mathbb{R})$. The conclusion follows from \eqref{bddeq1} and \eqref{bddeq2}.
	\endproof
	
	The illustration of Theorem \ref{bddseq} is provided in the following example.
	\begin{exa}\label{bddexa}
		Let $N=1$, consider the frames $\big\{\phi_{m,n}^{1,0}\big\}_{m, n\in\mathbb{Z}}$ for $L^2(\mathbb{R})$ and $\big\{\psi_{m,n}^{1,0}\big\}_{m, n\in\mathbb{Z}}$ given in Example \ref{finiteexa} and  Example \ref{operexa1}, respectively with frame bounds $A_1=4$, $B_1=16$ and $A_2= B_2=4$.
		
		Define the sequences $\{\alpha_{m,n}\}_{m,n\in \mathbb{Z}}$ and $\{\beta_{m,n}\}_{m,n\in \mathbb{Z}}$ as follows:
		\begin{align*}
			\alpha_{m,n} &= \frac{(-1)^{m n}}{2} \quad
			\text{and} \quad
			\beta_{m,n} = (-1)^n 3, \  \ m,n \in \mathbb{Z}.
		\end{align*}
		Therefore, $\inf |\alpha_{m,n}| = \sup  |\alpha_{m,n}| =\frac{1}{2}$ and $\inf |\beta_{m,n}| =\sup |\beta_{m,n}| = 3$, satisfying,
		\begin{align*}
			24 = 2 \sup |\alpha_{m,n}| \sup |\beta_{m,n}| \sqrt{B_1 B_2} < \inf |\alpha_{m,n}|^2 A_1 + \inf |\beta_{m,n}|^2 A_2= 37.
		\end{align*}
		Hence, by Theorem \ref{bddseq}, $\big\{\alpha_{m,n}\phi_{m,n}^{1,0}+ \beta_{m,n}\psi_{m,n}^{1,0} \big\}_{m, n \in \mathbb{Z}}$ is a frame for $L^2(\mathbb{R})$ with frame bounds $A= 13$ and $B=64$.
	\end{exa}

\section{Applications  to Frame Algorithm and Convergence Analysis}\label{subsec2}
In this section, we discuss applications of our results on sums of frames in the frame algorithm. Note that our results are also valid for discrete frames in separable Hilbert spaces. We recall that the frame decomposition of a frame gives stable reconstruction of each vector in the space. In the frame decomposition, we should compute the inverse of the frame operator, which may be complicated or take reasonable time if size of the matrix of the frame operator is large. In this situation, the frame algorithm is useful and used to approximate functions (signals). The frame algorithm  is based on the  width of the frame and hence on frame bounds of a given frame.
	The width of a frame $\mathcal{F}$ with frame bounds $A$,  $B$, denoted by $\Delta_{\mathcal{F}}$,  is  defined as
	\begin{align}\label{wid1}
		\Delta_{\mathcal{F}} := \text{width of} \ \mathcal{F} = \frac{B - A} {B + A}.
	\end{align}
	Thus, the width of the frame $\mathcal{F}$  measures its ``\emph{tightness}''. We also note that $\Delta < 1$.
	Nowadays there are a variety of  algorithms which depend on frame bounds and are used in approximations of signals  \cite{ole}.

\begin{thm}$[$Frame Algorithm$]$ \label{thm44}\cite[\  p. 14]{ole}
	Let $\mathcal{F} = \{\varphi_k\}_{k=1}^{m}$ be a frame for a finite dimensional  Hilbert space $\mathbb{H}$  with frame bounds  $A$, $B$ and frame operator $S$.  For a given  $\varphi \in \mathbb{H}$, define  a sequence $\{\psi_k\}_{k=0}^{\infty}$ in  $\mathbb{H}$ as follows.
	\begin{align*}
		\psi_0 =0, \quad \psi_k = \psi_{k-1} + \frac{2}{B + A} S(\varphi- \psi_{k-1}), \  \  k \in \mathbb{N}.
	\end{align*}
	Then, $\{\psi_k\}_{k=0}^{\infty}$ converges to $\varphi$ in $\mathbb{H}$ and the rate of convergence is given by
	\begin{align}\label{IIIeqthm3.1}
		\|\varphi - \psi_k\| \leq \big(\Delta_{\mathcal{F}}\big)^{k} \|\varphi\|, \ \  k \in \mathbb{N} \cup \{0\}, \  \text{where} \ \Delta_{\mathcal{F}} \  \text{is the width of} \ \mathcal{F}.
	\end{align}
\end{thm}

The following example gives an application of Theorem \ref{sfn} in the frame algorithm.
	
\begin{exa}\label{exfal}
 Consider the family of vectors  $\mathcal{F}: = \{f_k\}_{k=1}^{3}= \Big\{(\sqrt{6},\sqrt{6}),(0,2),(2,0)\Big\}$ in the complex space  $\mathbb{H}=\mathbb{C}^2$. Then,  a  simple computation shows that $\mathcal{F}$ is a frame for $\mathbb{H}$  with frame bounds $A_1=4$ and $B_1=16$. The width of $\mathcal{F}$ is given by
	\begin{align}\label{widI}
		\Delta_{\mathcal{F}}= \frac{B_1 - A_1}{B_1 + A_1} = \frac{3}{5}= 0.6.
	\end{align}
Now the family $\mathcal{F}_2:= \{g_k\}_{k=1}^{3} = \Big\{(0,3),(\sqrt{3},0),(\sqrt{3},0)\Big\} \subset \mathbb{H}$ is a tight frame for $\mathbb{H}$  with frame bounds $A_2 = B_2=9$.
	
	Choose $c_1=1$, $c_2=100$. Then, an estimate of condition  \eqref{sfc} in Theorem \ref{sfn} with $j=1$ is given by
	\begin{align*}
		 |c_1|A_1+\bigg|\frac{c_2^2}{c_1}\bigg|A_2 = 90004> 2400 = 2|c_2| \sqrt{B_1B_2}
	\end{align*}
	Hence, by Theorem \ref{sfn}, the sum $\mathcal{F}_{\sum}: = \left\{c_1f_k + c_2 g_k  \right\} _{k=1}^3$ is a frame for $\mathbb{C}^2$ with frame bounds \break
$A = \Big(\sum\limits_{i=1}^{2}|c_i|^2 A_i-2\sum\limits_{i=1 \atop i\neq j}^{2}|c_jc_i|\sqrt{B_j B_i}\Big) = 87604$ and $B = 2 \sum\limits_{i=1}^{2}|c_i|^2B_i = 180032$.
The width  of  the finite sum of frames $\mathcal{F}_{\sum}$  is given by
	\begin{align}\label{widIII}
		{\Delta_{\mathcal{F}}}_{\sum} = \frac{B - A}{B + A}=\frac{92428}{267636}=0.3453 \ (\text{up to} \ 4 \ \text{decimal places}).
	\end{align}
	From \eqref{widI} and \eqref{widIII}, we have ${\Delta_{\mathcal{F}}}_{\sum} < \Delta_{\mathcal{F}}$. Thus, the speed of convergence in the frame algorithm can be increased considerably if we use frame bounds of  finite linear combination of frames, see the graph below.\\
	\begin{figure}[H]
		\centering
		\includegraphics[width=3 in]{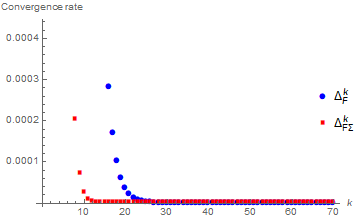}
		\caption{}
		\label{figfinite}
	\end{figure}
	
\end{exa}

The following example provides an application of Theorem \ref{dual} in the frame algorithm.
	\begin{exa}\label{dualexaalgo}
 Consider the frames  $\mathcal{F}:\{f_k\}_{k=1}^5=\Big\{\Big(\frac{1}{\sqrt3},0,0\Big),\Big(0,\frac{1}{\sqrt3},0\Big),\Big(0,0,\frac{1}{\sqrt3}\Big),\Big(\sqrt2,0,0\Big),\Big(0,0,\sqrt2\Big)\Big\}$  for the complex Hilbert space $\mathbb{H} = \mathbb{C}^3$ with frame bounds $A_1= \frac{1}{3}$ and $B_1= \frac{7}{3}$. Its dual frame is given by
 \begin{align*}
 \mathcal{G}:\{g_k\}_{k=1}^5=\Big\{\Big(\frac{\sqrt3}{2},0,0\Big),\Big(0,\sqrt3,0\Big),\Big(0,0,\frac{\sqrt3}{2}\Big),\Big(\frac{1}{2\sqrt2},0,0\Big),\Big(0,0,\frac{1}{2\sqrt2}\Big)\Big\}
 \end{align*}
with frame bounds $A_2 = \frac{7}{8}$ and $B_2 =3$. Using  Theorem \ref{dual}, the sum $\mathcal{F} + \mathcal{G}: =\big\{f_k+g_k\big\}_{k=1}^5$ is a frame with frame bounds $A= \frac{77}{24}$ and $B= \frac{22}{3}$.

The following table demonstrates that the width of $\mathcal{F} + \mathcal{G}$ is less than the  width of $\mathcal{F}$ and $\mathcal{G}$.
		
		\begin{center}
			\begin{tabular}{ |p{4.5cm}|p{2.5cm}|p{2.5cm}|p{2.5cm}|  }
				\hline
				\begin{center}
					\textbf{Frame}
				\end{center}& \begin{center}
					\textbf{Lower Bound}
				\end{center}&  \begin{center}
					\textbf{Upper Bound}
				\end{center}& \begin{center}
					\textbf{Width of Frame}  (up to $4$ decimal places)
				\end{center}\\
				\hline \begin{center}
					$\mathcal{F}: = \big\{f_k\big\}_{k=1}^5$\end{center} & \begin{center}$\frac{1}{3}$  \end{center} & \begin{center} $\frac{7}{3}$ \end{center}& \begin{center} $0.75$ \end{center}\\
				\hline
				\begin{center}
					$\mathcal{G}: = \big\{g_k\big\}_{k=1}^5$\end{center} & \begin{center}$\frac{7}{8}$  \end{center} & \begin{center} $3$ \end{center}& \begin{center} $0.5483$ \end{center}\\
				\hline \begin{center}
					$\mathcal{F} + \mathcal{G}: =\big\{f_k+g_k\big\}_{k=1}^5$\end{center} & \begin{center}$\frac{77}{24}$  \end{center} & \begin{center} $\frac{22}{3}$ \end{center}& \begin{center} $0.3913$ \end{center}\\
				\hline
			\end{tabular}
		\end{center}
The speed of convergence in the frame algorithm increases for the sum of a frame and its dual frame, see  Figure \ref{figdual} below.
		\begin{figure}[H]
			\centering
			\includegraphics[width=3 in]{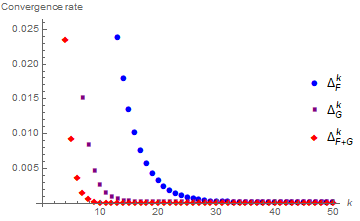}
			\caption{}
			\label{figdual}
		\end{figure}
	\end{exa}

 In the next example, we show that the frame bounds  of sum of images of two frames can increase the rate of approximation in the frame algorithm. This illustrates Theorem \ref{operator}.
	
	\begin{exa}\label{operexaalgo}
		Consider the frame $\mathcal{F}:=\{f_k\}_{k=1}^3=\{(\sqrt{6},\sqrt{6}),(0,2),(2,0)\}$ for $\mathbb{C}^2$ given in Example \ref{exfal} with frame bounds $A_1=4$ and $B_1=16$. Now the collection of vectors  $\mathcal{G}:=\{g_k\}_{k=1}^3=\{(2,0),(0,\sqrt{2}),(0,\sqrt{2})\}$   is a tight frame for $\mathbb{C}^2$ with frame bounds $A_2=B_2 =4$.
Define bounded linear operators $\Theta_1$, $\Theta_2 \colon \mathbb{C}^2 \rightarrow \mathbb{C}^2$ as follows:
		\begin{align*}
	\Theta_1 (z_1,z_2) =\Big(\frac{z_1}{160},\frac{z_2}{160}\Big) \quad \text{and} \quad
	\Theta_2 (z_1,z_2) =(z_1,z_2), \ (z_1,z_2) \in \mathbb{C}^2 .
		\end{align*}
		Then, it is easy to see that $\|\Theta_1^*\| =  \|\Theta_1\| = \frac{1}{160}$ and $\|\Theta_2^*\| =  \|\Theta_2\| = 1$. Further, $\|\Theta_1^*\|$ and $\|\Theta_2^*\|$ are bounded below by  $m_1 = \frac{1}{160}$ and $m_2 = 1$, respectively. Therefore
		\begin{align*}
			\frac{25601}{6400} = A_1m_1^2 + A_2 m_2^2 > 2\sqrt{B_1 B_2}\|\Theta_1\|\|\Theta_2\| = \frac{1}{10}.
		\end{align*}
		Hence, by Theorem \ref{operator}, $\Theta_1 \mathcal{F}+\Theta_2 \mathcal{G}: =\{\Theta_1f_k+\Theta_2g_k\}_{k=1}^3$ is a frame for $\mathbb{C}^2$ with frame bounds $A = \frac{24961}{6400}$ and $B = \frac{6561}{1600}$.
		
The width of the frames $\mathcal{F}$ and  $\Theta_1 \mathcal{F}+\Theta_2 \mathcal{G}$ are given below.
		\begin{align*}
\Delta_{\mathcal{F}} = \frac{B_1-A_1}{B_1+A_1} = \frac{3}{5}= 0.6 \quad \text{and} \quad \Delta_{\Theta_1 \mathcal{F}+\Theta_2 \mathcal{G}}= \frac{B-A}{B+A} =\frac{1283}{51205}= 0.0250 \ (\text{up to} \ 4 \ \text{decimal places}).
		\end{align*}
This shows that the width of  the frame $\Theta_1 \mathcal{F}+\Theta_2 \mathcal{G}$	is considerably less than   the  width  of the frame $\mathcal{F}$.
	\end{exa}

We conclude the section with an application of Theorem \ref{bddseq}.
	\begin{exa}\label{bddexaalgo}
		Consider the frames  $\mathcal{F}:=\{f_k\}_{k=1}^3$ and $\mathcal{G}:=\{g_k\}_{k=1}^3$ for the space $\mathbb{C}^2$ given in  Example \ref{operexaalgo} with frame bounds $A_1=4$, $B_1=16$ and $A_2= B_2=4$, respectively. Let $\alpha: = \{\alpha_k\}_{k=1}^3=\Big\{\frac{(-1)^{k}}{4}\Big\}_{k=1}^3$ and $\beta:=\{\beta_k\}_{k=1}^3=\Big\{(-1)^k 4\Big\}_{k=1}^3$. Then, $\min\limits_{1 \le k \le 3} |\alpha_k| =\frac{1}{4}=\max\limits_{1 \le k \le 3} |\alpha_k| $ and $\min\limits_{1 \le k \le 3} |\beta_k| = 4=\max\limits_{1 \le k \le 3}|\beta_k| $, satisfying
		\begin{align*}
			16 = 2 \max\limits_{1 \le k \le 3} |\alpha_k| \max\limits_{1 \le k \le 3} |\beta_k| \sqrt{B_1 B_2} < \min\limits_{1 \le k \le 3} |\alpha_k|^2 A_1 + \min\limits_{1 \le k \le 3} |\beta_k|^2 A_2= \frac{257}{4}.
		\end{align*}
		Hence, by Theorem \ref{bddseq}, the weighted sum $\alpha \mathcal{F} + \beta \mathcal{G}: = \big\{\alpha_k f_k+ \beta_k g_k \big\}_{k=1}^3$ is a frame for $\mathbb{C}^2$ with frame bounds $A=\frac{193}{4}$ and $B=81$.

The width of the frames $F$ and $\alpha \mathcal{F} + \beta \mathcal{G}$ are given in the following equations:
		\begin{align*}
			\Delta_{\mathcal{F}}= \frac{B_1-A_1}{B_1+A_1} = \frac{3}{5}= 0.6 \quad \text{and} \quad
			\Delta_{\alpha \mathcal{F} + \beta \mathcal{G}} = \frac{B-A}{B+A} =\frac{131}{517}= 0.2533 \ (\text{up to} \ 4 \ \text{decimal places}).
		\end{align*}
One may observe  the difference between the width of  $\mathcal{F}$ and $\alpha \mathcal{F} + \beta \mathcal{G}$, and hence the rate of convergence in the frame algorithm using frame bounds of $\mathcal{F}$ and the weighted sum of frames $\alpha \mathcal{F} + \beta \mathcal{G}$.	
	\end{exa}

\section{Discussion and Conclusion}
Some of the  most significant applications of frames involve approximation of signals. Frame algorithms depend on the width of the frame which measures tightness of the given frame and provides the rate of approximation of signals. The width of a frame is a function of frame bounds of a given frame. Using frame bounds of finite sums of frames, the width of the frame can be decreased and hence the rate of convergence in the frame algorithms increases.  For this reason, we study sums of frames  of the signal space $L^2(\mathbb{R})$ which are associated with the Weyl-Heisenberg group. By using idea given in  \cite{JDLG}  for sums of frames with wave packet structure in  matrix-valued function spaces, our result  Theorem \ref{sfn} provides new  sufficient conditions for finite sums of frames from the Weyl-Heisenberg group to be a frame of the space $L^2(\mathbb{R})$. In \cite{OSCT}, the authors proved necessary and sufficient  conditions, in terms of invertibility of a certain type of operator, for sum of images of two Bessel sequences; and sum of image a frame and its dual   under bounded linear operators on a Hilbert space to be a frame for the underlying Hilbert space.  Without using invertibility of operators associated with given frames and Bessel sequences, in  Theorem \ref{dual},  we  provide explicit frame bounds of sum of a  frame and its dual frame from the Weyl-Heisenberg group. Theorem \ref{operator} gives new sufficient conditions, with explicit frame bounds,  under which sum of images of two frames for $L^2(\mathbb{R})$  under bounded linear operators on $L^2(\mathbb{R})$  forms a frame for  $L^2(\mathbb{R})$. In Theorem \ref{bddseq}, we  present frame conditions for  perturbed sums of frames, that is, sums of frames   where frame vectors are multiplied, or perturbed by  bounded sequences of scalars.

Majority of research on sums of frames focuses on necessary and sufficient conditions under which finite sum of frames turns out to be a frame for the underlying space. Our finding, in the application part of this paper, shows that  frame bounds of  sums of frames can decrease the width of the frame which increases the rate of approximation of signals in the frame algorithm. The details are  given in  Example \ref{exfal}, Example \ref{dualexaalgo}, Example \ref{operexaalgo} and Example \ref{bddexaalgo} which illustrate Theorem \ref{sfn}, Theorem \ref{dual}, Theorem \ref{operator} and Theorem \ref{bddseq}, respectively. Further, the rate of convergence is illustrated in Fig.1 and Fig.2. We believe that our work will be useful in applications of frames in approximation of signals from the frame algorithm. The present study still leaves many directions for applications of frame bounds of sums of frames  in other variants of frame algorithms, for example Chebyshev and conjugate gradient method \cite{KGalgo}.

%

\end{document}